\documentclass[envcountsame,envcountsect]{svmult}

\usepackage{amssymb}
\usepackage{amsmath}

\newtheorem{notation}{Notation}

\newtheorem{examples}{Examples}

\begin{document}

\title*{Heat Kernels and Critical Limits}
\titlerunning{Heat Kernels}

\author{Doug Pickrell}
\institute{Mathematics Department, University of Arizona \\
Tucson, AZ, USA 85721 \\
\texttt{pickrell@math.arizona.edu}}

\maketitle

\begin{abstract} This article is an exposition of several questions
linking heat kernel measures on infinite dimensional Lie groups,
limits associated with critical Sobolev exponents, and
Feynmann-Kac measures for sigma models. The first part of the
article concerns existence and invariance issues for heat kernels.
The main examples are heat kernels on groups of the form
$C^0(X,\mathbf F)$, where $X$ is a Riemannian manifold and
$\mathbf F$ is a finite dimensional Lie group. These measures
depend on a smoothness parameter $s>dim(X)/2$. The second part of
the article concerns the limit $s\downarrow dim(X)/2$, especially
$dim(X)\le 2$, and how this limit is related to issues arising in
quantum field theory. In the case of $X=S^1$, we conjecture that
heat kernels converge to measures which arise naturally from the
Kac-Moody-Segal point of view on loop groups, as $s\downarrow
1/2$.
\end{abstract}

\section{Introduction}
\label{sec0}

Given a finite dimensional real Hilbert space $\mathbf f$, there
is an associated Lebesgue measure $\lambda_{\mathbf f}$, a
(positive) Laplace operator $\Delta_{\mathbf f}$, and a
convolution semigroup of heat kernel measures
\begin{equation}\nu_t^{\mathbf f}=e^{-\frac t2\Delta_{\mathbf f}}\delta_0=(2\pi t)^{
-dim(\mathbf f)/2}e^{-\frac 1{2t}\vert x\vert_{\mathbf
f}^2}d\lambda_{ \mathbf f}(x).\end{equation} The map $\mathbf f\to
\{\nu_t^{\mathbf f}\}$ is functorial, in the sense that if
$P:\mathbf f_1\to \mathbf f_2$ is an orthogonal projection, then
\begin{equation}P_{*}\nu_t^{\mathbf f_1}=\nu_t^{\mathbf f_2}.\label{0.2}\end{equation}

More generally, given a finite dimensional Lie group $\mathbf F$,
with a fixed inner product on its Lie algebra $\mathbf f$, there
is an induced left invariant Riemannian metric, a normalized Haar
measure $\lambda_{\mathbf F}$, a Laplace operator $
\Delta_{\mathbf F}$, and a convolution semigroup of heat kernel
measures,
\begin{equation}\nu_t^{\mathbf F}=e^{-t\Delta_{\mathbf F}/2}\delta_1=\lim_{N\to\infty}
(exp_{*}(\nu_{t/N}^{\mathbf f}))^N\end{equation} (the power refers
to convolution of measures on $\mathbf F$). The map $\mathbf F\to
\{\nu_t^{\mathbf F}\}$ is functorial, in the sense that if
$P:\mathbf F_1\to \mathbf F_2$ and the induced map of Lie algebras
$\mathbf f_ 1\to \mathbf f_2$ is an orthogonal projection, then
\begin{equation}P_{*}\nu_t^{\mathbf F_1}=\nu_t^{\mathbf F_2}.\end{equation}

In Section 2 we will recall how this construction can be
generalized to infinite dimensions.  The form of the linear
generalization is well-known:  given a separable real Hilbert
space $\mathbf h$, there is an associated convolution semigroup of
Gaussian measures $\{\nu_t^{\mathbf h}\}$.  An essential
complication arises when $dim(\mathbf h)=\infty$, namely the
$\nu_t^{ \mathbf h}$ are only finitely additive, viewed as
cylinder measures on $\mathbf h$ (see \cite{St} for orientation on
this point).  To overcome this, following Gross, we consider an
abstract Wiener space $\mathbf h\to \mathbf b$; in this framework,
the Gaussians are realized as countably additive measures
$\nu_t^{\mathbf h\subset \mathbf b}$ on $ \mathbf b$, a Banach
space completion of $\mathbf h$. Geometrically speaking, one
imagines that the support of the heat kernel diffuses into the
enveloping space $\mathbf b$.

In the larger group category, the objects are abstract Wiener
groups:  such an object consists of an inclusion of a
separable Hilbert Lie group into a Banach Lie group,
\begin{equation}\mathbf H\to \mathbf G,\end{equation}
such that the induced map of Lie algebras $\mathbf h\to \mathbf g$
is an abstract Wiener space. With some possible restrictions, the
corresponding heat kernels $\nu_t^{\mathbf H\subset \mathbf G}$
can be constructed as in the finite dimensional case, and they
form a convolution semigroup of inversion-invariant probability
measures on $\mathbf G$ (this use of Ito's ideas in an infinite
dimensional context apparently originated in \cite{DS}). This
construction is functorial, in the same sense as in finite
dimensions.

In Section 3 we consider invariance properties of heat kernels.
The Cameron-Martin-Segal theorem asserts that for $t>0$, the
measure class $[\nu_t^{\mathbf h\subset \mathbf b}]$ is invariant
with respect to translation by $h\in \mathbf b$ if and only if
$h\in \mathbf h$.  In the larger group category, we conjecture
that $[\nu_t^{\mathbf H\subset \mathbf G}]$ is invariant with
respect to translation by $h\in \mathbf G$ if and only if $h\in
\mathbf H_0$ and $Ad(h)^{*}Ad(h)-1$ is a Hilbert-Schmidt operator
on $\mathbf h$. The evidence in favor of this conjecture is based
on deep results of Driver.

In Section 4, and in the remainder of the paper, we focus on
groups of maps. Suppose that $X$ is a compact manifold and $W$ is
a real Hilbert Sobolev space corresponding to degree of smoothness
$s>dim(X)/2$. In this context the Sobolev embedding
\begin{equation}W\to C^0(X;\mathbb R)\end{equation}
is a canonical example of an abstract Wiener space. If $\mathbf F$
is a finite dimensional Lie group with a fixed inner product on
its Lie algebra $\mathbf f$, then
\begin{equation}\mathbf H=W(X,\mathbf F)\to \mathbf G=C^0(X,\mathbf F)\label{0.7}\end{equation}
is an abstract Wiener group.  In this example the measure
$\nu_t^{\mathbf H\subset \mathbf G}$ is determined by its
distributions corresponding to evaluation at a finite number of
points of $X$; these distributions involve a nonlocal Green's
function interaction between pairs of points.

When $X=S^1$, $\mathbf F=K$, a simply connected compact group with
$Ad(K)$-invariant inner product, and $s=1$, Driver and
collaborators have proven that the measure class of
$\nu_t^{\mathbf H\subset \mathbf G}$ equals the measure class of a
Wiener measure. This latter measure has a heuristic expression in
terms of a local kinetic energy functional
\begin{equation}exp(-\frac 1{2t}\int_X\langle g^{-1}dg\wedge *g^{-1}dg\rangle )
\prod_Xd\lambda (g(\theta ))\label{0.8}\end{equation} (see
\cite{Dr2}). When $dim(X)=2$, $s=1$ is the critical exponent, the
heat kernels are defined for $s>1$, and the heuristic expression
(\ref{0.8}) is the Feynmann-Kac measure for the quantum field
theory with fields $g:X\to K$ (the sigma model with target $K$).

This motivates the following questions, considered in Sections 5
and 6: is there an analogue of (\ref{0.7}) when $s=dim(X)/2$, and
is there a sense in which heat kernels have limits as $s\downarrow
dim(X)/2$? When $s=dim(X)/2$, the Sobolev embedding fails (because
of scale invariance), and fields cannot be evaluated at individual
points. Moreover for fields with values in a curved space such as
$K$, there is not an obvious notion of generalized function.

In the case of $X=S^1$, there is at least one possible resolution.
In addition to accepting that the support diffuses outside of
ordinary functions, one also has to accept that it diffuses in
noncompact directions in $G$, the complexification of $K$. For
various reasons (e.g. the Kac-Moody-Segal point of view on loop
groups) one is led to consider, in place of (\ref{0.7}), the
hyperfunction completion
\begin{equation}W^{1/2}(S^1,K)\to Hyp(S^1,G).\end{equation}
A generic point in the hyperfunction completion is
represented by a formal Riemann-Hilbert factorization
\begin{equation}g=g_{-}\cdot g_0\cdot g_{+},\label{10}\end{equation}
where $g_{-}\in H^0(\Delta^{*},0;G,1)$, $g_0\in G$, $g_{+}\in H^0
(\Delta ,0;G,1)$, and $\Delta$ ($\Delta^{*}$) is the open unit
disk centered at $0$ ($\infty$).  There is a family of measures on
$Hyp(S^1,G)$, indexed by a level $l$, having heuristic expressions
for their densities involving Toeplitz determinants. Whereas heat
kernels (and Wiener measures) are determined by distributions
corresponding to points on $S^1$, these latter measures are in
theory determined by distributions corresponding to points off of
$S^1$; in reality, an enlightening way to display these measures
has not been found. Based on symmetry considerations, we
conjecture that the heat kernels parameterized by $s$ and $t$
converge to a measure with level $l=1/2t$ as $s\downarrow 1/2$.

The case of most interest is $X=\hat{\Sigma}$, a closed Riemannian
surface, because of (\ref{0.8}) and the central role of sigma
models in Physics. This is intertwined with the one dimensional
case, by the trace map $W^1(\hat{\Sigma })\to W^{1/2} (S)$, when
$S$ is a $1$-manifold (space) embedded in $\hat{\Sigma}$ (a
Euclidean space-time). In section 6 I have tried to point out a
few mathematical facts about sigma models which might be relevant
to understanding how to formulate a critical limit as $s\downarrow
1$.

Surveys of related topics, with more extensive bibliographies,
include \cite{Dr2} (for Wiener measure and heat kernels on loop
groups), \cite{A} (for energy representations), and \cite{Ga} (for
sigma models).

\begin{notation}(1) Throughout this paper, all spaces are assumed to be
separable and real, unless explicitly noted otherwise.

(2) We will frequently encounter continuous maps of Hilbert spaces
$\mathbf h\to \mathbf f$ such that the quotient Hilbert space
equals the Hilbert space structure on the target.  If $\mathbf f$
is a subspace of $\mathbf h$ itself, then such a map is an
orthogonal projection.  In this paper we will refer to such a map
$\mathbf h\to \mathbf f$ as a projection.

(3) Suppose that $\mathbf h$ is a Hilbert space.  Let $\mathcal
L_2$ denote the symmetrically normed ideal of Hilbert-Schmidt
operators on $\mathbf h$.  For an invertible operator $A$ on
$\mathbf h$, the two conditions
\begin{equation}AA^{-t}\in 1+\mathcal L_2, \quad
and \quad A-A^{-t}\in \mathcal L_2\end{equation} are equivalent.
We let $GL(\mathbf h)_{(\mathcal L_2)}$ denote the group of
invertible operators satisfying these two conditions (in the same
way, we can define $GL(\mathbf h)_{(\mathcal I)}$, for any
symmetrically normed ideal $\mathcal I$).

(4) Given a measure $\nu$ on a space $X$, and a Borel space $Y$,
a map $f:X\to Y$ will be said to be $\nu$-measureable if $f^{-1}(
E)$
is $\nu$-measureable for each Borel set $E\subset Y$.  The measure
class of $\nu$ will be denoted by $[\nu ]$.

(5) Lebesgue measure will be denoted by $d\lambda$, and Haar
measure will be denoted by $d\lambda_G$, for a group $G$.

\end{notation}

\section{General Constructions}
\label{sec1}

\subsection{Abstract Wiener spaces and Gaussian
measures} \label{sec1.1}

Suppose that $\mathbf h$ is a Hilbert space.  The corresponding
convolution semi-group of Gaussian measures, $\{\nu^{\mathbf
h}_t\}_{ t>0}$, has the heuristic expression
\begin{equation}d\nu_t^{\mathbf h}=(2\pi t)^{-dim(\mathbf h)/2}e^{-\frac 1{2t}\vert
x\vert^2_{\mathbf h}}d\lambda_{\mathbf
h}(x),\label{1.12}\end{equation} where $\lambda_{\mathbf h}$
denotes the `Lebesgue measure for the Hilbert space $\mathbf h$'.
The measure $\lambda_{\mathbf h}$ and the expression (\ref{1.12})
have literal meaning only when $\mathbf h$ is finite dimensional.

To go beyond this, following Gross, suppose that $\mathbf h\to
\mathbf b$ is a continuous dense inclusion of $\mathbf h$ into a
Banach space. We can consistently define a finitely additive
probability measure $\nu_t^{\mathbf h\subset \mathbf b}$ on
$\mathbf b$ in the following way.

A cylinder set is a Borel subset of $\mathbf b$ of the form
$p^{-1} (E)$, where $p:\mathbf b\to \mathbf f$, $p\vert_{\mathbf
h}:\mathbf h\to \mathbf f$ is a finite rank projection, and $ E$
is a Borel subset of $\mathbf f$.  The set of cylinder sets is an
algebra; it is a $\sigma$-algebra if and only if $\mathbf
h=\mathbf b$ is finite dimensional.  The functorial property of
finite dimensional Gaussian measures, (\ref{0.2}), implies that
there exists a well-defined finitely additive measure
$\nu_t^{\mathbf h\subset \mathbf b}$ on the algebra of cylinder
sets satisfying
\begin{equation}p_{*}\nu_t^{\mathbf h\subset \mathbf b}=\nu_t^{\mathbf f},\end{equation}
for all $p$ as above.

\begin{definition} (a) The inclusion $\mathbf h\to \mathbf b$ is called an
abstract Wiener space if the finitely additive measure
$\nu_t^{\mathbf h\subset \mathbf b}$ has a (necessarily unique)
countably additive extension to a Borel probability measure on
$\mathbf b$, for some $t>0$ (and hence all $t>0$, because
$\nu_t^{\mathbf h\subset \mathbf b} (E)=\nu_1^{\mathbf h\subset
\mathbf b}(E/\sqrt t)$).

(b) A map of abstract Wiener spaces,
\begin{equation}(\mathbf h_1\subset \mathbf b_1)\to (\mathbf h_2\subset \mathbf b_2)\end{equation}
is a map $\mathbf b_1\to \mathbf b_2$ such that the restriction to
$\mathbf h_ 1$ is a projection $\mathbf h_1\to \mathbf h_2$.

(c) The cylindrical Schwarz space of $\mathbf h\subset \mathbf b$
is
\begin{equation}\mathcal S(\mathbf h\subset \mathbf b)=\lim_pp^{*}\mathcal S(\mathbf f),\end{equation}
where the limit is a directed limit over all finite rank maps of
Wiener spaces $p:\mathbf b\to \mathbf f$.

(d) The Laplace operator acting on a cylindrical Schwarz
function is given by
\begin{equation}\Delta_{\mathbf h\subset \mathbf b}(p^{*}\phi )=p^{*}(\Delta_{\mathbf f}
\phi )\end{equation}
\end{definition}

For a map of abstract Wiener spaces, the corresponding
Gaussian semigroups push forward.

In many ways, the restriction to Banach space completions of
$\mathbf h$ is artificial. However the Banach framework seems
elegant and natural, in part because of the following
characterization and examples.

\begin{theorem} The inclusion $\mathbf h\to \mathbf b$ is an abstract
Wiener space if and only if for each $\epsilon >0$ there exists a
finite dimensional subspace $\mathbf f_{\epsilon}\subset \mathbf
h$ such that
\begin{equation}\nu_t^{\mathbf h\subset \mathbf b}\{x\in \mathbf b:\vert px\vert_{\mathbf f}
>\epsilon \}<\epsilon\end{equation}
for all $p:\mathbf b\to \mathbf f$ vanishing on $\mathbf
f_{\epsilon}$, where $ p:\mathbf h\to \mathbf f$ is a finite rank
projection.  If $\mathbf b$ is itself a Hilbert space, this is
equivalent to the condition that the inclusion $\mathbf h\to
\mathbf b$ is a Hilbert-Schmidt operator.
\end{theorem}

For a discussion of this theorem, see section 3.9 of \cite{B}.

\begin{examples}
 (a) If $\mathbf f$ is a finite dimensional Hilbert space, then
$\mathbf f=\mathbf h=\mathbf b$ is an abstract Wiener space, and
$\nu_t^{\mathbf f}$ is given literally by (\ref{1.12}).

(b) Suppose that $X$ is a compact manifold, possibly with
boundary, and $W$ is a Hilbert Sobolev space of degree of
smoothness $s>dim(X)/2$.  Then the Sobolev embedding
\begin{equation}W\to C^0(X;\mathbb R)\end{equation}
is an abstract Wiener space.  This is closely related to
the fact that $W\to W^0$ is a Hilbert-Schmidt operator
precisely when
$s$ is above the critical exponent.

(c) For $W$ as in (b), if $\mathbf h\to \mathbf b$ is an abstract
Wiener space, then
\begin{equation}W(X,\mathbf h)=W\otimes \mathbf h\to C^0(X,\mathbf b)\end{equation}
is also an abstract Wiener space (see $3.11.29$ of \cite{B}).

\end{examples}

I do not know of a good reference for (b). However the basic ideas
are in \cite{CL}. A special case of (c) is the Wiener space
analogue of the path functor,
\begin{equation}W\otimes \mathbf h=W^1([0,T],0;\mathbf h,0)\to C^0([0,T],0;\mathbf b,
0)\label{20}\end{equation} with inner product
\begin{equation}\langle x,y\rangle_{W^1}=\langle\dot {x},\dot {y}\rangle_{L^2}=
\int_0^T\langle\dot {x}(\tau ),\dot {y}(\tau )\rangle_{\mathbf
h}d\tau .\label{21}\end{equation} To stay within the Banach
category, we have choosen $T<\infty$.  By letting $T\to\infty$, we
obtain a Gaussian semigroup on the semi-infinite path space
\begin{equation}C^0([0,\infty ),0;\mathbf b,0).\end{equation}
For notational simplicity, we will write
\begin{equation}\nu^{\mathbf h\subset \mathbf b}=\nu_1^{W\otimes \mathbf h\subset C^0
([0,\infty ),0;\mathbf b,0)},\end{equation} and we will refer to
this as the Brownian motion associated to $\mathbf h\subset
\mathbf b$.  Given a partition $V:0<t_ 1<..<t_n$, there is an
evaluation map
\begin{equation}Eval_V:C^0([0,\infty ),0;\mathbf b,0)\to\prod_1^n\mathbf b:x\to (x_
i),\quad x_i=x(t_i).\end{equation} In terms of the coordinates
$x_i$,
\begin{equation}(Eval_V)_{*}\nu^{\mathbf h\subset \mathbf b}=d\nu^{\mathbf h\subset \mathbf b}_{
t_1}(x_1)\times d\nu^{\mathbf h\subset \mathbf
b}_{t_1-t_2}(x_2-x_1)\times ..\times d\nu^{\mathbf h\subset
\mathbf b}_{t_n-t_{n-1}}(x_n-x_{n-1}) .\end{equation} In
particular the distribution of $\nu^{\mathbf h\subset \mathbf b}$
at time $ t$ is $\nu_t^{\mathbf h\subset \mathbf b}$. The
semigroup property of $\nu_t^{\mathbf h\subset \mathbf b}$ is
equivalent to the consistency of these distributions for
$\nu^{\mathbf h\subset \mathbf b}$ when the partition is refined.

The path space construction can be iterated. Thus given $\mathbf
h\subset \mathbf b$, there is a Brownian motion, a Brownian sheet,
etc., associated to $\mathbf h\subset \mathbf b$.

\subsection{Abstract Wiener groups and heat
kernels} \label{1.2}

\begin{definition} (a) An abstract Wiener group is an
inclusion of a separable Hilbert Lie group $\mathbf H$ into a
Banach Lie group $\mathbf G$ such that the associated Lie algebra
map $\mathbf h\to \mathbf g$ is an abstract Wiener space.

(b) A map of abstract Wiener groups is a map of Lie groups
$\mathbf G_1\to \mathbf G_2$ such $\mathbf H_1\to \mathbf H_2$ and
the induced Lie algebra map is a projection.

(c) An abstract Wiener group $\mathbf H\subset \mathbf G$ with the
property that finite rank maps of Wiener groups $p:\mathbf G\to
\mathbf F$ separate the topology of $\mathbf G$ is called local.

(d) For a local abstract Wiener group, we define the
cylindrical Schwarz space by
\begin{equation}\mathcal S(\mathbf H\subset \mathbf G)=\lim_pp^{*}\mathcal S(\mathbf F)\end{equation}
and the (left-invariant) Laplacian acting on a cylindrical
Schwarz function by
\begin{equation}\Delta_{\mathbf H\subset \mathbf G}(p^{*}\phi )=p^{*}\Delta_{\mathbf F}
\phi .\end{equation}
\end{definition}

Given an abstract Wiener group $\mathbf H\subset \mathbf G$, there
is an induced left-invariant Riemannian metric on $\mathbf H$, and
a left-invariant Finsler metric on $\mathbf G$.  These induce
separable complete metrics compatible with the topologies of
$\mathbf H$ and $\mathbf G$.

\begin{examples}(a) If $\mathbf F$ is a finite dimensional Lie group, with an
arbitrary inner product on its Lie algebra $\mathbf f$, then
$\mathbf F=\mathbf H=\mathbf G$ is an abstract Wiener group.

(b) If $W$ is as in (b) of Examples 1, given an abstract Wiener
group $\mathbf H\subset \mathbf G$,
\begin{equation}W(X,\mathbf H)\subset C^0(X,\mathbf G)\end{equation}
is a local abstract Wiener group.

(c) The construction in (b) can be adapted to local gauge
transformations of a nontrivial principal bundle $P\to X$,
provided that the structure group $K$ is a compact Lie
group with a fixed $Ad(K)$-invariant inner product.

(d) There are also nonlocal examples involving infinite classical
matrix groups; see \cite{Go}.

\end{examples}

The following would be the ideal existence result.

\begin{conjecture}For fixed $t\ge 0$, the sequence of
probability measures
\begin{equation}exp_{*}(\nu_{t/N}^{\mathbf h\subset \mathbf g})^N\label{29}\end{equation}
(the $N$-fold convolution on $\mathbf G$) has a weak limit with
respect to $BC(\mathbf G)$, bounded continuous functions. The
limits, denoted $\nu_t^{\mathbf H\subset \mathbf G}$, form a
convolution semigroup of inversion-invariant probability measures
on $\mathbf G$.  Given a map
\begin{equation}P:(\mathbf H_1\subset \mathbf G_1)\to (\mathbf H_2\subset \mathbf G_2)\end{equation}
of abstract Wiener groups, $P_{*}\nu_t^{\mathbf H_1\subset \mathbf
G_ 1}=\nu_t^{\mathbf H_2\subset \mathbf G_2}$, for each $t\ge 0$.
\end{conjecture}

This conjecture is known to be true if one additionally assumes
that the inclusion $\mathbf h\subset\mathbf g$ is 2-summable (see
\cite{BD}).

I have not stated this result in a standard form, and this
obscures an important idea. To explain this idea, fix $t>0$.
Consider the map
\begin{equation}C^0([0,\infty ),0;\mathbf g,0)\to \mathbf G:x\to g_N(t)=e^{x_1}e^{x_2-
x_1}..e^{x_N-x_{N-1}},\end{equation} where $x_i=x(it/N)$.  The
$g_N(t)$-distribution of the Brownian motion measure $\nu^{\mathbf
h\subset \mathbf g}$ is the $N$-fold convolution (\ref{29}),
because the $x_i-x_{i-1}$ are independent ($\mathbf g$-valued)
Gaussian variables with variance $t/N$.  Thus the conjecture is
equivalent to the assertion that $g_N(t)$ has a limit as $N \to
\infty$, say $g(t)$, with $\nu^{\mathbf h\subset \mathbf
g}$-probability one. The latter statement implies that the heat
kernel is an image of a Gaussian measure, an insight due to Ito
(the standard approach goes a step further, and constructs $g(t)$,
as a continuous function of $t$, by showing that it solves a
stochastic differential equation).

\begin{examples}  (a) If $\mathbf F$ is a finite
dimensional real Lie group, with an arbitrary inner product on the
Lie algebra $\mathbf f$, then for $\Delta$ the Laplacian for the
left invariant Riemannian metric,
\begin{equation}\nu_t^{\mathbf F}=e^{-t\Delta /2}\delta_1,\end{equation}
This measure is absolutely continuous with respect to
Haar measure, and has a real analytic density.

(b) Consider the path construction
\begin{equation}\mathbf H=W^1(I,0;\mathbf F,1)\to \mathbf G=C^0(I,0;\mathbf F,1)\end{equation}
where $I=[0,T]$, and $\mathbf h$ has the norm (\ref{21}) (with
$\mathbf f$ in place of $\mathbf h$ in (\ref{20})). In this case
$\nu_t^{\mathbf H\subset \mathbf G}$ is Brownian motion for $
\mathbf F$, with variance $t$. The path space for $\mathbf G$ is
the double path space
\begin{equation}C^0(I\times I,\{0\}\times I\cup I\times \{0\};\mathbf F,1),\end{equation}
and $\nu^{\mathbf H\subset \mathbf G}$ is the $\mathbf F$-valued
Brownian sheet.

\end{examples}

The Brownian sheet is relevant to Physics in the following way.
Supppose that $\mathbf F=K$, a compact Lie group, and the inner
product is $Ad(K)$-invariant. Suppose that $\Sigma$ is a closed
surface with an area element. The Yang-Mills measure is the
measure on the space of $K$-connections which has the heuristic
expression
\begin{equation}e^{-\frac 12\int\langle F_A\wedge *F_A\rangle}d\lambda (A)\end{equation}
for $A\in\Omega^1(\Sigma ;\mathfrak k)$, where $F_A$ is the
curvature for the connection $d+A$, and $d\lambda$ denotes a
heuristic Lebesgue measure on the linear space of $\mathfrak
k$-valued one forms. This measure can be (roughly) defined as a
finitely additive measure in the following way. Let
$p_t(g)d\lambda (g)$ denote the heat kernel for $K$. Given a
triangulation of $\Sigma$, one considers the projection
\begin{equation}\Omega^1(\Sigma ,\mathfrak k)\to K^E:A\to (g_e)\end{equation}
where $E$ denotes the set of edges for the triangulation,
and $g_e\in K$ represents parallel translation for the
connection $d+A$. The image of the Yang-Mills measure is
\begin{equation}\prod_Fp_{tArea(f)}(g_{\partial f})\prod_Edg_e,\end{equation}
where $g_{\partial f}$ denotes the holonomy around the face $f\in
F$. The fact that these measures are coherent with respect to
refinement of the triangulation follows directly from the
semigroup property of the heat kernel. When projected to gauge
equivalence classes, this measure can be expressed as a finite
dimensional conditioning of the Brownian sheet (see \cite{Sen}).

\section{Invariance Questions}
\label{sec2}

The linear invariance properties of Gaussian measures are
summarized by the following

\begin{theorem} Suppose that $\mathbf h\to \mathbf b$ is an abstract
Wiener space.  Let $\nu =\nu_t^{\mathbf h\subset \mathbf b}$ for
some $ t>0$.

(a) The subset of translations in $\mathbf b$ which fix $[\nu ]$
is $ \mathbf h$.

(b) The group of linear transformations which fix $[\nu ]$ is
$GL(\mathbf h)_{(\mathcal L_2)}$, in the following sense:  given a
$\nu$-measureable linear map $\hat {T}:\mathbf b\to \mathbf b$
which fixes $[\nu ]$, $T =\hat {T}\vert_{\mathbf h}\in GL(\mathbf
h)_{(\mathcal L_2)}$; conversely, given $T\in GL(\mathbf
h)_{(\mathcal L_2)}$, there exists a $\nu$-measureable linear map
$\hat {T}:\mathbf b\to \mathbf b$ which fixes $ [\nu ]$ and
satisfies $\hat {T}\vert_{\mathbf h}=T$.

\end{theorem}

These results, and their history, are discussed in a very
illuminating way in \cite{B} (where one can also find expressions
for the Radon-Nikodym derivatives).

Now suppose that $\mathbf H\subset \mathbf G$ is an abstract
Wiener group, and $g\in \mathbf G$.  When is
\begin{equation}[(L_g)_{*}\nu_t^{\mathbf H\subset \mathbf G}]=[\nu_t^{\mathbf H\subset
\mathbf G}],\label{38}\end{equation} where $L_g$ is left
translation? Necessarily $g\in \mathbf G_0$, the identity
component (because $\nu_ t^{\mathbf H\subset \mathbf G}$ is
supported in the identity component).  The abelian case suggests
that $g\in \mathbf H$. Because $\nu_t^{\mathbf H\subset \mathbf
G}$ is inversion invariant, translation invariance implies that
$[\nu_t^{\mathbf H\subset \mathbf G}]$ is $g$-conjugation
invariant.  Thus if (\ref{38}) holds for all $t$, this suggests
that $[\nu_t^{\mathbf h\subset \mathbf g}]$ is $Ad(g)$-invariant,
and this implies that, as an operator on $\mathbf h$, $Ad(g)\in
GL(\mathbf h)_{(\mathcal L_2)}$.  These considerations suggest the
following

\begin{conjecture} For $T>0$ and $h\in \mathbf G$, $[\nu^{\mathbf H
\subset \mathbf G}_T]$ is left and right $h$-invariant if and only
if $h\in \mathbf H_0$ and $Ad(h)\in GL(\mathbf h)_{(\mathcal
L_2)}$.
\end{conjecture}

Hilbert-Schmidt criteria are universal in
unitary representation-theoretic questions of this type.
In the next section we will consider evidence
in favor of this conjecture.

Before leaving this abstract setting, I will mention
two other related questions.

When $t\to\infty$, a Gaussian measure is asymptotically invariant
in the following sense introduced by the Malliavins: given
$\nu_t=\nu_t^{\mathbf h\subset \mathbf b}$ and $ h\in \mathbf h$,
for each $p<\infty$,
\begin{equation}\int\vert 1-\frac
{d\nu_t(b+h)}{d\nu_t(b)}\vert^pd\nu_t(b)\to 0\end{equation} as
$t\to\infty$ (in this linear context, this integral reduces to a
one-dimensional integral, which is easily estimated). Now suppose
that Conjecture 3.2 is true, and let $\nu_t=\nu_t^{\mathbf
H\subset \mathbf G}$ and $h$ satisfy the conditions in Conjecture
3.2. Is it true that for each $p<\infty$,
\begin{equation}\int\vert 1-\frac {d\nu_t(gh)}{d\nu_t(g)}\vert^p d\nu_t(g)\to 0\end{equation}
as $t\to\infty$? The Malliavins proved that Wiener measure on
$C^0(S^1,K)$ is asymptotically invariant in this sense (see
chapter 4, Part III, of \cite{Pi1} for a quantitative version of
this result, and an example of how this result is used to prove
existence of an invariant measure for the loop group).

The second question is the following. Suppose that we change the
inner product on $\mathbf h$ in a way which (by (b) of Theorem
3.1) does not change the measure class of $\nu_t^{\mathbf h
\subset \mathbf g}$:
\begin{equation}\langle x,y\rangle_1=\langle Ax,y\rangle ,\quad x,y\in \mathbf h\end{equation}
where $A-1\in \mathcal L_2(\mathbf h)$. Write $\mathbf h_1$ for
$\mathbf h$ with this inner product depending on $A$.  Is
$[\nu_t^{\mathbf H_1\subset \mathbf G}] =$ $[\nu^{\mathbf H\subset
\mathbf G}_t]$, for each $t$?

\section{The Example of $Map(X,\mathbf F)$}
\label{sec3}

Let $X$ denote a compact Riemannian manifold, and fix a
Sobolev space (with a fixed inner product) of continuous
real-valued functions
\begin{equation}W\hookrightarrow C^0(X),\end{equation}
necessarily corresponding to some degree of smoothness
$s>dim(X)/2$. For the sake of clarity, we will
suppose that $\partial X$ is empty, and the inner product
for $W$ is of the form
\begin{equation}\langle\phi_1,\phi_2\rangle_W=\int_X(P^{s/2}(\phi_1))(P^{s/2}(\phi_
2))dV\end{equation} where $P$ denotes a positive elliptic
psuedodifferential operator of order $2$, e.g.  $P=(m_0^2+\Delta
)$, where $\Delta$ is a (positive) Laplace type operator.  We will
write $G(x,y)$ for the Green's function of $P^s$; thus
\begin{equation}P^s_y(G(x,y)dV(y))=\delta_x(y)\end{equation}
in the sense of distributions. Evaluation at $x\in X$ defines a
continuous linear functional
\begin{equation}\delta_x:W\to \mathbb R:\phi\to\phi (x);\end{equation}
this functional is represented by $G(x,\cdot )\in W$:
\begin{equation}\phi (x)=\langle\phi ,G(x,\cdot )\rangle_W.\label{3.46}\end{equation}

Let $\mathbf F$ denote a finite dimensional abstract Wiener group.
Then
\begin{equation}\mathbf H=W(X,\mathbf F)\to \mathbf G=C^0(X,\mathbf F)\label{47}\end{equation}
is an abstract Wiener group, where $\mathbf h=W\otimes \mathbf f$,
as a Hilbert space.  It is local because we can evaluate at points
of $X$.

Given a finite set of points $V\subset X$, let
\begin{equation}Eval_V:C^0(X,\mathbf F)\to \mathbf F^V:g\to (g(v))\end{equation}
denote the evaluation map.  This is a map of abstract Wiener
groups, where $\mathbf f^V$ has the Hilbert space structure
induced by the Lie algebra map
\begin{equation}Eval_V:W(X,\mathbf f)\to \mathbf f^V.\end{equation}
It follows from (\ref{3.46}) that this inner product on $\mathbf
f^V$ is given by
\begin{equation}\langle (x_v),(y_w)\rangle =\sum_{v,w\in V}G^{v,w}\langle x_v,y_
w\rangle_{\mathbf f},\label{50}\end{equation} where $(G^{v,w})$ is
the matrix inverse to the covariance matrix $(G(v,w))_{v,w\in V}$.

More generally,
given an embedded submanifold $S\subset X$,
restriction induces a map of abstract Wiener groups
\begin{equation}C^0(X,\mathbf F)\to C^0(S,\mathbf F),\end{equation}
where the inner product on $W(S,\mathbf f)$ corresponds to a
degree of smoothness $s-codim(S)/2$.  This is very reminiscent of
functorial properties which are exploited at a heuristic level for
Feynmann-Kac measures (see Section 6).

If $f_1,..,f_N$ denotes an orthonormal basis for $\mathbf f$, and
if $f_i^{(v)}$ denotes the left invariant vector field on $\mathbf
F^V$ corresponding to the $v$-coordinate, then the left invariant
Laplacian on $\mathbf F^V$ determined by the inner product
(\ref{50}) is given by
\begin{equation}\Delta_V=\sum_{v,w\in V}\sum_{j=1}^NG(v,w)f_i^{(v)}f_j^{(w)},\end{equation}
and
\begin{equation}(Eval_V)_{*}\nu_t^{\mathbf H\subset \mathbf G}=e^{-t\Delta_V/2}\delta_
1^{(\mathbf F^V)}.\label{53}\end{equation} These finite
dimensional distributions determine the measure $\nu_t^{\mathbf
H\subset \mathbf G}$. We will compare this with Wiener measure
below, when $X=S^1$.

We now want to understand the content of Conjecture 3.2. Suppose
that $g\in W(X,\mathbf F)$.  We need to compute the adjoint for
$Ad(g)$, denoted $Ad(g)^{*_W}$, acting on the Hilbert Lie algebra
$W(X,\mathbf f)$.  Let $Ad(g)^t$ denote the adjoint for $Ad(g)$
acting on $L^2(X,dV)\otimes \mathbf f$.  Then

\begin{equation}\langle Ad(g)^{*_W}(f_1),f_2\rangle_{\mathbf h}=\int P^sf_1Ad(g)(
f_2)dV\end{equation}
\begin{equation}=\int P^s(P^{-s}Ad(g)^tP^s)(f_1)f_2dV.\end{equation}
Thus
\begin{equation}Ad(g)^{*_W}=P^{-s}\circ Ad(g)^t\circ P^s,\end{equation}
and
\begin{equation}Ad(g)+Ad(g)^{*_W}=(Ad(g)P^{-s}+P^{-s}Ad(g)^t)P^s.\end{equation}
In general this is simply a zeroth order operator. However suppose
that $Ad(g)$ has values in $O(\mathbf f)$, the orthogonal group.
Then
\begin{equation}Ad(g)+Ad(g)^{*_W}=[Ad(g),P^{-s}]P^s,\label{58}\end{equation}
and because of the commutator, the order drops to $-1$ (see
\cite{Fr}).

To formulate this in operator language, let $\mathcal L_d^{+}$
denote the ideal of compact operators $T$ on $\mathbf h$
satisfying
\begin{equation}\sup_N\frac 1{log(N)}\sum_1^Ns_n(T),\quad\end{equation}
where $s_1\ge s_2\ge ..$ are the eigenvalues of $\vert T\vert$.
The fact that (\ref{58}) has order $-1$ implies that it belongs to
$\mathcal L_d^{+}$. This implies the following statement.

\begin{theorem} Let $d=dim(X)$.

(a) Suppose that $g\in W(X,\mathbf f)$ and $Ad(g(x))\in O(\mathbf
f)$, for each $x\in X$.  Then
\begin{equation}Ad(g)\in GL(\mathbf h)_{(\mathcal L_d^{+})}.\end{equation}
(b) Suppose that $\mathbf F=K$, or $\mathbf F=K^{\mathbb C}$,
where $K$ is a compact Lie group, and the inner product on
$\mathbf f$ is $Ad(K)$-invariant.  Assuming the truth of
Conjecture 3.2, $[\nu_t^{\mathbf H\subset \mathbf G}]$ is
$W(X,K)$-invariant when $ X=S^1$, and noninvariant for $dim(X)>1$.
\end{theorem}

The important point is that $dim(X)=2$ is the critical
case, where the criterion marginally fails.

We now want to consider evidence in support of Conjecture 3.2.

Suppose that $X=S^1$, and $P=m_0^2-(\frac {\partial}{\partial\theta}
)^2$. For $s>1/2$,
\begin{equation}G(\theta_1,\theta_2)=G(\theta_1-\theta_2)=\frac 1{\pi}\sum_{n\ge
0}(m_0^2+n^2)^{-s}cos(n(\theta_1-\theta_2)).\end{equation}

We will write $\nu_t^{s,m_0}$ for the corresponding heat kernels.
When $s=1$, we
can compare $\nu_t^{s,m_0}$ with Wiener measure $\omega_t$. This latter
measure has distributions
\begin{equation}(Eval_V)_{*}\omega_t=\prod_Ep_t(g_vg_w^{-1})\prod_Vd\lambda_K(g_
v),\end{equation} where $E$ denotes the set of edges (between
vertices in $V$), and $\nu_t^K=p_td\lambda_K$. The important point
is that this is a nearest neighbor interaction, which is far more
elementary than the interaction involving all pairs of points for
heat kernels, as in (\ref{53}).

The following summarizes some of the deep
results of Driver and collaborators on
heat kernels for loop groups.

\begin{theorem}Suppose that $X=S^1$ and $\mathbf F=K$,
a simply connected compact Lie group with $Ad(K)$-invariant inner
product on $\mathbf f$.

(a) For $s=1$, $[\nu_t^{s,m_0}]=[\omega_t]$.

(b) For all $s>1/2$, $[\nu_t^{s,m_0}]$ is $W(S^1,K)$-biinvariant.
\end{theorem}

A relatively short and illuminating discussion of (a) can be found
in \cite{Dr2}. Part (b) is a long story. For $s=1$, part (b)
follows from (a) and the relatively well-understood fact that
$[\omega_t]$ is $W^1(S^1,K)$-biinvariant. Driver gave a direct
proof of (b), for $s=1$ (\cite{Dr1}), and the ideas were shown to
extend to all $s>1/2$ in \cite{I}.

In support of Conjecture 3.2, one can also directly analyze the
Brownian motion construction ((b) of Example 2), using a
stochastic Ito map (which essentially linearizes the problem).

When $X$ is a closed Riemannian surface, we expect that the heat
kernels for $W^{s,m_0}(X,K)\to C^0(X,K)$ to just miss being
translation quasiinvariant, for all $s>1$. It is interesting to
ask whether there is another natural construction which yields a
quasiinvariant measure. This should be compared with Shavgulidze's
construction of quasiinvariant measures for $Diff(X)$, for $X$ of
any dimension (see page 312 of \cite{B}).

\section{The Critical Sobolev Exponent and $X=S^1$}
\label{sec4}

For $s=dim(X)/2$, the Sobolev embedding fails, and in general a
$W^{dim(X)/2}$-function on $X$ is not bounded (for references to
the finer properties of generic functions with critical Sobolev
smoothness, see \cite{H}). This leads to a subtle situation. On
the one hand, while $W^{dim(X)/2}(X,\mathfrak k)$ is a Hilbert
space, unless $\mathfrak k$ is abelian, it is not a Lie algebra.
On the other hand, because $K$ is assumed compact, a map $g\in
W^{dim(X)/2}(X,K)$ is automatically bounded. Consequently
$W^{dim(X)/2}(X,K)$, with the Sobolev topology (we do not impose
uniform convergence), is a topological group, but it is not a Lie
group (the study of the algebraic topology of spaces with
topologies defined by critical Sobolev norms is nascent, see
\cite{Br}).

For nonabelian $\mathfrak k$, in general it is simply not clear
how to form an analogue of (\ref{47}), when $s=dim(X)/2$.

In the rest of this section we will focus on the case $X=S^1$. We
will freely use facts about loop groups, as presented in
\cite{PS}. For simplicity of exposition, we will assume that $K$
is simply connected, $\mathfrak k$ is simple, and the inner
product (on the dual) satisfies $\langle\theta ,\theta \rangle=2$,
where $\theta$ is a long root (the abelian case is essentially
trivial, but requires qualifications). The complexification of $K$
will be denoted by $G$.

For each $s>1/2$, there is a (Kac-Moody) universal central
extension of Lie groups
\begin{equation}0\to \mathbb T\to\hat {W}^s(S^1,K)\to W^s(S^1,K)\to 0.\label{63}\end{equation}
This extension exists for $s=1/2$ as well (as an extension of a
topological group), and this is the maximal domain for the
extension. This extension is realized using operator methods in
\cite{PS}. If we fix a faithful representation $K\subset U(N)$, a
loop $g:S^1\to K$ can be viewed as a unitary multiplication
operator $M_g$ on $H=L^2(S^1,\mathbb C^N)$. Relative to the Hardy
polarization $H=H_{+}\oplus H_{-}$, where $H_{+}$ consists of
functions with holomorphic extension to the disk,
\begin{equation}M_g=\left(\begin{array}{cc} A&B\\
C&D\end{array} \right),\end{equation} where $B$ (or $C$) is the
Hankel operator, and $A$ (or $D$) is the Toeplitz operator,
associated to the loop $g$. The condition $g\in W^{1/2}$ is
equivalent to $B\in \mathcal L_2$. The Toeplitz operator defines a
map
\begin{equation}A:W^{1/2}(S^1,K)\to Fred(H_{+}).\end{equation}
The lift of the left action of $W^{1/2}(S^1,K)$ on itself to
$A^{*}Det,$ the pullback of the determinant line bundle $Det\to
Fred(H_{+})$, induces (a power of) the extension (\ref{63}). The
upshot is that, from an analytic point of view, Kac-Moody theory
is intimately related to the critical exponent and Toeplitz
determinants.

The group $H^0(S^1,G)$ is
a complex Lie group.  An open, dense neighborhood of the
identity consists of those loops which have a unique
(triangular or Birkhoff or Riemann-Hilbert) factorization
\begin{equation}g=g_{-}\cdot g_0\cdot g_{+},\label{66}\end{equation}
where $g_{-}\in H^0(D^{*},\infty ;G,1)$, $g_0\in G$, $g_{+}\in H^
0(D,0;G,1)$, and $D$ and $D^{*}$ denote the $\underline {closed}$
unit disks centered at $ 0$ and $\infty$, respectively.  A model
for this neighborhood is
\begin{equation}H^1(D^{*},\mathfrak g)\times
G\times H^1(D,\mathfrak g),\end{equation} where the linear
coordinates are determined by $\theta_{+}=g_{+}^{-1}(\partial
g_{+})$, $\theta_{-}=(\partial g_{ -})g_{-}^{-1}$.  The (left or
right) translates of this neighborhood cover $H^0(S^1,G)$.

The hyperfunction completion, $Hyp(S^1,G)$, is modelled on the
space
\begin{equation}H^1(\Delta^{*},\mathfrak g)\times G\times H^1(\Delta ,\mathfrak g),\end{equation}
where $\Delta$ and $\Delta^{*}$ denote the $\underline {open}$
disks centered at $ 0$ and $\infty$, respectively, and the
transition functions are obtained by continuously extending the
transition functions for the analytic loop space.  The group
$H^0(S^1,G)$ acts naturally from both the left and right of
$Hyp(S^1,G)$ (see ch. 2, Part III of \cite{Pi1}).

There is a holomorphic line bundle $\mathcal L\to Hyp(S^1,G)$ with
a map
\begin{equation}\begin{array}{ccc} \hat {W}^{1/2}(S^1,K)&\to&\mathcal L\\
\downarrow&&\downarrow\\
W^{1/2}(S^1,K)&\to&Hyp(S^1,G)\end{array} \end{equation} which is
equivariant with respect to natural left and right actions by
$\hat {H}^0(S^1,K)$, where $H^0(S^1,K)$ denotes real analytic maps
into $K$. The line bundle $\mathcal L^{*}$ has a distinguished
holomorphic section $\sigma_0$; restricted to $W^{1/2}$ loops,
$\mathcal L^{*}{}^ m=A^{*}Det$ and $\sigma_0^m=detA$, the pullback
of the canonical holomorphic section of $Det\to Fred(H_{+})$,
where $m$ is the ratio of the $\mathbb C^N$-trace form and the
normalized inner product on $\mathfrak k$. We will write $\vert
\mathcal L\vert^2$ for the line bundle $\mathcal L\otimes\bar
{\mathcal L}$; we can form real powers of this bundle, because the
transition functions are positive.

\begin{theorem}For each $l\ge 0$, there exists
a $H^0(S^1,K)$-biinvariant measure $d\mu^{\vert \mathcal
L\vert^{2l}}$ with values in the line bundle $\vert \mathcal
L\vert^{2l}\to Hyp(S^1,G)$ such that
\begin{equation}d\mu_l=(\sigma_0\otimes\bar{\sigma}_0)^ld\mu^{\vert \mathcal L\vert^{
2l}}\end{equation} is a probability measure. In particular there
exists a $H^0(S^1,K)$-biinvariant probability measure $d\mu
=d\mu^{\vert \mathcal L\vert^0}$ on $Hyp(S^ 1,G)$.
\end{theorem}

This is a refinement of the main result in \cite{Pi1} (see also
\cite{Pi2}). These bundle-valued measures are conjecturally
unique. Uniqueness would imply invariance with respect to the
natural action of analytic reparameterizations of $S^1$. This
independence of scale is the hallmark of the critical exponent.

One motivation for constructing the measure $\mu^{\vert \mathcal
L\vert^{2l}}$ is to prove a Peter-Weyl theorem of the schematic
form
\begin{equation}H^0\cap L^2(\mathcal L^{*}{}^{\otimes l})=\sum_{level(\Lambda )=l}H
(\Lambda )\otimes H(\Lambda )^{*}\label{71}\end{equation} where
the $H(\Lambda )$ are the positive energy representations of level
$l\in \mathbb Z^{+}$ (This statement requires more explanation.
Here we will simply say, Kac and Peterson proved an algebraic
generalization of the Peter-Weyl theorem (see \S 1.7, Part I of
\cite{Pi1}), and we are seeking an analytic version, suitable for
application to sewing rules discussed below). This remains
incomplete.

Because $\vert\sigma_0\vert^2=det\vert A\vert^{2/m}$, restricted to $
W^{1/2}(S^1,K)$,
$\mu_l$ has a heuristic expression
\begin{equation}d\mu_l=\frac 1{\mathcal Z}det\vert A\vert^{2(l/m)}d\mu .\label{72}\end{equation}
This begs two questions: (1) why do Toeplitz
determinants have anything to do with (limits of) heat kernels
and (2) how do we write the background $\mu$ in a
way which suggests how to think about
it analytically?

It is instructive to first consider the abelian case. When
$K=\mathbb T$, the theorem is valid provided $l>0$ (the background
$\mu$ does not exist, reflecting a lack of `compactness' in this
flat case), and we consider identity components. An element of
$Hyp(S^1,\mathbb C^{*})_0$ can be written uniquely as
\begin{equation}exp(\sum_{n<0}x_nz^n)\cdot exp(x_0)\cdot exp(\sum_{n>0}x_nz^n)\label{73}\end{equation}
where the sums represent holomorphic functions in the open disks.
If (\ref{73}) represents a loop $g\in W^{1/2}(S^1,\mathbb T)_0$,
then
\begin{equation}det\vert A(g)\vert^{2l}=exp(-l\sum_{n>0}n\vert x_n\vert^2)\end{equation}
(the Helton-Howe formula). Consequently
\begin{equation}d\mu_l=d\lambda_{\mathbb T}(exp(x_0))\prod_{k>0}\frac 1{\mathcal Z_k}e
xp(-lk\vert x_k\vert^2)d\lambda (x_k),\label{75}\end{equation} and
there is no diffusion in noncompact directions: $x_{-k}=-x_k^{*}$
$a.e.$. Since
\begin{equation}d\nu_t^{s,m_0}=d\nu_{m_0^{-2s}t}^{\mathbb T}(exp(x_0))\prod_{k>0}\frac
1{\mathcal Z_k}exp(-\frac 1{2t}(k^2+m_0^2)^s\vert
x_k\vert^2)d\lambda (x_k),\end{equation} (\ref{75}) is the limit
of heat kernels (with $l=1/2t$) as $s\downarrow 1 /2$, provided
that we also send the mass $m_0$ to zero.

Now suppose that $K$ is simply connected.
Since $K$ and its loop group fit into the
Kac-Moody framework, it is natural
to compare $\mu$ to Haar measure $d\lambda_K$ for $K$.

Fix a triangular decomposition
\begin{equation}\mathfrak g=\mathfrak n^{-}\oplus \mathfrak h\oplus \mathfrak n^{+}.\end{equation}
A generic $g\in K$ (or $G$), can be written uniquely as
\begin{equation}g=lmau,\end{equation}
where $l\in N^{-}=exp(\mathfrak n^{-})$, $m\in T=exp(\mathfrak
h\cap \mathfrak k )$, $a\in A=exp(\mathfrak h_{\mathbb R})$, and
$u\in N^{+}=exp(\mathfrak n^{+})$.  This implies that there is a
$K$-equivariant isomorphism of homogeneous spaces, $K/T \to
G/B^{+}$ and $l\in N^{-}$ is a coordinate for the top stratum.
Harish-Chandra discovered that in terms of this coordinate, the
unique $K$-invariant probability measure on $K/T$ can be written
as
\begin{equation}a^{4\delta}d\lambda (l)=\prod\vert\sigma_i(g)\vert^2d\lambda (l
)=\frac 1{\vert l\cdot v_{\delta}\vert^4}d\lambda
(l),\label{79}\end{equation} where (\ref{78}) implicitly
determines $a=a(gT)$ as a function of $l$, the $\sigma_i$ are the
fundamental matrix coefficients, $d\lambda (l)$ denotes a properly
normalized Haar measure for $N^{-}$, $2\delta$ denotes the sum of
the positive complex roots for the triangular factorization
(\ref{77}), and $v_{\delta}$ is a highest weight vector in the
highest weight representation corresponding to the dominant
integral functional $\delta$.  The point of the third expression
is that it shows the denominator for the density is a polynomial
in $l$, hence the integrability of (\ref{79}) is quite sensitive.
It also follows from work of Harish-Chandra that for $\lambda\in
\mathfrak h_{\mathbb R}^{*}$
\begin{equation}\int a^{-i\lambda}d\lambda_Kg=\prod_{\alpha >0}\frac {\langle 2
\delta ,\alpha\rangle}{\langle 2\delta -i\lambda
,\alpha\rangle}\label{80}\end{equation} (the right hand side is
Harish-Chandra's $\mathbf c$-function for $G/K$). In particular
this determines the integrability of powers of $a$: for $\epsilon
>0$,
\begin{equation}\int a^{(1+\epsilon )2\delta}d\lambda (l)=\epsilon^{-(d-r)/2}<\infty
,\end{equation} $d$ and $r$ denote the dimension and rank of
$\mathfrak g$, respectively.  The formula (\ref{79}) can be
derived by calculating a Jacobian in a straightforward way. The
formula (\ref{80}) can be deduced from the Duistermaat-Heckman
localization principle (in particular $log(a)$ is a momentum map),
using a (Drinfeld-Evens-Lu) Poisson structure which generalizes to
the Kac-Moody framework (see \cite{Pi3}).

Now consider the loop situation. Recall that $\theta_{-}$ is a
coordinate for $H^0(\Delta^{*},0;G,1)$ (which is similar to
$N^{-}$). The loop analogue of the formula (\ref{79}) leads to the
following heuristic expression for the $\theta_{-}$ distribution
of $\mu_l$, where initially we think of $\theta_{-}$ as
corresponding to a unitary loop $g\in W^{1/2}(S^1,K)$:
\begin{equation}(\theta_{-})_{*}\mu_l=\frac 1{\mathcal Z}\vert\sigma_0(g)\vert^{2(2
\dot {g}+l)}d\lambda (\theta_{-})\end{equation}
\begin{equation}=\frac 1{\mathcal Z}det\vert A(g)\vert^{2(2\dot {g}+l)/m}d\lambda (
\theta_{-})\label{83}\end{equation} where $\dot {g}$ is the dual
Coxeter number, and $d\lambda (\theta_{-})$ is a heuristic
Lebesgue background. The important point is the shift by the dual
Coxeter number, which reflects a regularization of the `sum over
all positive roots'. In a similar manner the Duistermaat-Heckman
approach to (\ref{80}) applies in a heuristic way to compute the
$g_0$ distribution of $\mu_l$:
\begin{equation}\int a(g_0)^{-i\lambda}d\mu_l=\prod_{\alpha >0}\frac {sin(\frac
1{2(\dot {g}+l)}\langle 2\delta ,\alpha\rangle )}{sin(\frac
1{2(\dot { g}+l)}\langle 2\delta -i\lambda ,\alpha\rangle
)}.\label{84}\end{equation} The function on the right hand side is
an affine analogue of Harish-Chandra's $\mathbf c$-function.

The heuristic expression (\ref{83}) explains why $\mu_l$ is
expected to be invariant with respect to $PSU(1,1)$ (the Toeplitz
determinant and the background Lebesgue measure are conformally
invariant). It also points to a way of constructing the
$\theta_{-}$ (or $g_{-}$) distribution of $\mu_l$, by imposing a
cutoff and taking a limit:
\begin{equation}d((g_{-})_{*}\mu_l)(\theta_{-})=\lim_{N\uparrow\infty}\frac 1{\mathcal Z_{
P_N}}det\vert A(g(P_N\theta_{-}))\vert^{2(2\dot {g}+l)/m}d\lambda
(P_N\theta_{-}),\end{equation} where $P_N\theta_{-}$ is the
projection onto the first $N$ coefficients, $g_{-}$ is related to
$P\theta_{-}$ by $P\theta_{-} =(\partial g_{-})g_{-}^{-1}$, and
$g$ is a unitary loop having Riemann-Hilbert factorization
(\ref{66}). Many, but not all, of the details of this construction
have been worked out. The main idea is that $log(det\vert
A(g)\vert^2)$ is part of a momentum map. Consequently localization
can be used to compute integrals involving the determinant against
a symplectic volume element, on finite dimensional approximations.

The formula (\ref{84}) for the $g_0$ distribution remains a
conjecture, but it seems to point in a fruitful direction. The
zero-mode $g_0$ roughly arises from the projection
\begin{equation}H^0(\hat {\mathbb C},0,\infty ;G,1)\to Hyp(S^1,G)_{generic}\to G\end{equation}
More generally,
given a closed Riemann surface $\hat{\Sigma}$ and a
real analytic embedding $S^1\to\hat{\Sigma}$, a $G$-hyperfunction induces
a holomorphic $G$-bundle on $\hat{\Sigma}$, and consequently there is
a bundle
\begin{equation}H^0(\hat{\Sigma}\setminus S^1,G)\to Hyp(S^1,G)\to H^1(\hat{\Sigma }
,\mathcal O_G),\label{86}\end{equation} where $H^1(\hat{\Sigma
},\mathcal O_G)$ denotes the set of isomorphism classes of
holomorphic $G$-bundles. The stable points of $H^1(\hat{\Sigma
},\mathcal O_G)$ identify with the irreducible points of
$H^1(\hat{\Sigma },K)$, which has a canonical symplectic structure
(see \cite{AB}). This suggests the following

\begin{question} Does $d\mu$ project to the (normalized) symplectic
volume on the stable points of $H^1(\hat{\Sigma },\mathcal O_G)$?
\end{question}

There is a natural generalization of this to include $d\mu^{\vert
\mathcal L\vert^{2l}}$. But it is not clear how to incorporate
reduction at points, as in (\ref{86}).

Continuous loops have Riemann-Hilbert factorizations,
and as a consequence there is a diagram
\begin{equation}\begin{array}{ccc} W^s(S^1,K)&\to&C^0(S^1,K)\\
\downarrow&&\downarrow\\
W^{1/2}(S^1,K)&\to&Hyp(S^1,G)\end{array} \end{equation} for each
$s>1/2$. We consider the $W^s$ norm
\begin{equation}\sum_n(m_0^2+n^2)^{s/2}\vert x_n\vert^2,\end{equation}
and $\nu_t^{s,m_0}$, the corresponding heat kernels, which we
view as measures on $Hyp(S^1,G)$.

\begin{conjecture} The measures $\nu_t^{s,m_0}$ have a limit as
$s\downarrow 1/2$, and the measure class of this limit equals the
measure class of $d\mu_l$, where $l=1/2t$. When
$m_0\downarrow 0$, this limit converges to $\mu_l$.
\end{conjecture}

The intuition is simply that for
$m_0=0$, the heat kernel should relax to a configuration
which is conformally invariant, as $s\downarrow 1/2$. The
existing hard evidence is slight. It is true in the
abelian case. The heat kernel measure class $[\nu_t^{s,m_0}]$
is biinvariant with respect
to $W^s(S^1,K)$, $s>1/2$. It is expected that $[\mu_l]$ will be biinvariant
with respect to $W^{1/2}(S^1,K)$, in the sense that the
natural induced unitary representation
\begin{equation}H^0(S^1,K\times K)\to U(L^2(d\mu_l))\end{equation}
will extend continuously to $W^{1/2}(S^1,K)$ (this is true for the
discrete part of the spectrum in (\ref{71})).

\section{2D Quantum Field Theory}
\label{sec5}

In two dimensions the question of how to formulate a notion of
critical limit for heat kernels is possibly related to the
question of how to mathematically formulate sigma models. To give
the discussion some structure, I will focus on one aspect of
quantum field theory (qft), following Segal.

Fix a dimension $d$, thought of as the dimension of (Euclidean)
space-time.  As in section 4 of \cite{Se}, let $\mathcal
C_{metric}$ denote the category for which the objects are oriented
closed Riemannian $(d-1)$-manifolds, and the morphisms are
oriented compact Riemannian $d$-manifolds with totally geodesic
boundaries.

\begin{definition} A primitive $d$-dimensional
unitary qft is a representation of $\mathcal C_{metric}$ by
separable Hilbert spaces and Hilbert-Schmidt operators such that
disjoint union corresponds to tensor product, orientation reversal
corresponds to adjoint, and $\mathcal C_{metric}$-isomorphisms
correspond to natural Hilbert space isomorphisms.
\end{definition}

It is an interesting question to what extent this definition
captures the meaning of locality in qft. Segal has recently
advocated additional axioms, and a priori there could be primitive
theories which are not truly local.

To my knowledge there does not exist an example of a nonfree qft
which has been shown to satisfy Segal's primitive axioms in
dimension $d>2$. However, in two dimensions, the space of all
theories is definitely large. From one point of view, this space
is the configuration space of string theory, and it is expected to
have a remarkable geometric structure; see \cite{Ma}.

Before discussing sigma models, we will recall, in outline, how a
`generic' theory is constructed, using constructive qft techniques
(see \cite{Pi4}).

Fix (a bare mass) $m_0>0$, and a polynomial $P:\mathbb R\to
\mathbb R$ which is bounded below. Given a closed Riemannian
surface $\hat{\Sigma}$, the $P(\phi )_2$-action is the local
functional
\begin{equation}\mathcal A:\mathcal F(\hat{\Sigma })\to \mathbb R:\phi\to\int_{\hat{\Sigma}}
(\frac 12(\vert d\phi\vert^2+m_0^2\phi^2)+P(\phi
))dA,\end{equation} where $\mathcal F(\hat{\Sigma })$ is the
appropriate domain of $\mathbb R$-valued fields on $\hat{\Sigma}$
for $\mathcal A$. A heuristic expression for the $P(\phi
)_2$-Feynmann-Kac measure is
\begin{equation}exp(-\mathcal A(\phi ))\prod_{x\in\hat{\Sigma}}d\lambda (\phi (x)).\label{92}\end{equation}
In this two dimensional setting, there is a rigorous
interpretation of (\ref{92}) as a finite measure on generalized
functions,
\begin{equation}\frac 1{det_{\zeta}(m_0^2+\Delta )^{1/2}}e^{-\int_{\hat{\Sigma}}
:P(\phi ):_{C_0}}d\phi_C,\label{93}\end{equation} where
$C=(m_0^2+\Delta )^{-1}$, $d\phi_C$ is the Gaussian probability
measure with covariance $C$ (corresponding to the critical Sobolev
space $W^1(\hat{\Sigma })$), $\int :P(\phi ):_{C_0}$ denotes a
regularization of the nonlinear interaction, and $det_{\zeta}$
denotes the zeta function determinant.

Suppose that $S$ is a closed Riemannian $1$-manifold. Given an
inner product for $W^{1/2}(S)$, there is an associated Gaussian
measure on generalized functions. We only need the measure class,
and hence we only require that the principal symbol of the
operator defining the inner product is compatible with length on
$S$. Associated to this measure class, say $\mathcal C(S)$, there
is a Hilbert space $\mathcal H( S)$, the space of half-densities
associated to $\mathcal C(S)$ (see the appendix to \cite{Pi4}).

Suppose that $\Sigma$ is a Riemannian surface with totally
geodesic boundary $S$. We consider the closed Riemannian surface
$\hat{\Sigma }=\Sigma^{*}\circ\Sigma$ obtained by gluing $ \Sigma$
to its mirror image $\Sigma^{*}$ along $S$. Because $d\phi_C$
corresponds to the critical exponent $s=1$, typical configurations
for the Feynmann measure (\ref{93}) are not ordinary functions.
However typical configurations are sufficiently regular so that it
makes sense to restrict them to $S$. Consequently the Feynmann
measure can be pushed forward to a finite measure on generalized
functions along $S$. Essentially because the underlying map of
Hilbert spaces is the trace map $W^1(\hat{\Sigma })\to
W^{1/2}(S)$, this finite measure lands in the measure class
$\mathcal C(S)$. By taking a square root (which is required
because $\Sigma$ is half of $\hat{\Sigma}$), we obtain a half
density $\mathcal Z(\Sigma )\in \mathcal H(S)$.

\begin{theorem}The maps $S\to \mathcal H(S)$ and $\Sigma\to \mathcal Z
(\Sigma )$ define
a representation of Segal's category.
\end{theorem}

This construction extends to vector space-valued fields
in a routine way. But severe problems arise for
qfts with classical fields having values in
nonlinear targets.

Consider the sigma model with target $K$. A naive idea is to use
heat kernels, in place of Gaussians, as backgrounds for an
analogous construction. This introduces a new element: we must
nudge $s>1$, then take a critical limit. From this point of view,
the renormalization group (see section 3 of \cite{Ga}) should be a
prescription for how to take this limit.

We will augment the sigma model action to include a topological
term, the Wess-Zumino-Witten `B-field'. Ignoring
$s$-regularization temporarily, the above outline reads as
follows: we represent the Feynmann-Kac measure
\begin{equation}exp(-\frac 1{2t}\int_{\hat{\Sigma}}\langle g^{-1}dg\wedge *g^{-
1}dg\rangle +2\pi
ilWZW(g))\prod_xd\lambda_K(g(x))\label{94}\end{equation} as a
density against a critical heat kernel background; we use the
naturality of heat kernels to push the Feynmann-Kac measure
forward along a map of generalized abstract Wiener spaces, induced
by the trace map $W^1(\hat{\Sigma },K)\to W^{1/2}(S,K)$,
\begin{equation}(W^1(\hat{\Sigma },K)\subset ?)\to (W^{1/2}(S,K)\subset Hyp(S,G
));\end{equation} and we obtain a Hilbert space of half-densities,
and a vector, by taking the positive square root of the
pushforward measure.

The point of introducing the $WZW$ term is that when the level $l$
is a positive integer, and $t=1/2l$, this ($WZW_l$) model is
solvable in many senses (see chapter 4 of \cite{Ga}, and
\cite{We}). It is believed that the model satisfies Segal's axioms
(see the Foreword to \cite{Se}). There also is the belief, far
more speculative, that there should be a one parameter family of
theories which interpolates between $WZW_l$ and the original sigma
model. This deformation intuitively arises by letting $t\downarrow
0$ in the action, and is related to the expected emergence of a
scale parameter in the process of regularizing the action
(\ref{94}), using for example heat kernels for $s>1$, and then
taking a limit as $s\downarrow 1$.

There are strong points of contact between our outline and what is
known about the $WZW_l$ model: The Hilbert space of the $WZW_l$
model is given by the right hand side of (\ref{71}). The vacuum of
the theory (the vector corresponding to a disk) is $det(A)^l$,
which is a "holomorphic square root" of $d\mu_l$ (this is far more
complicated than the positive square root for $P(\phi )_2$). The
proof of sewing (in the holomorphic sectors, involving the $WZW_l$
modular functor, see page 468 of \cite{Se}), ultimately hinges on
a generalized Peter-Weyl theorem (our conjectural analytical
version (\ref{71}) would fit perfectly with Segal's global
approach).

The bare sigma model is believed also to be solvable, but in a
completely different sense: in terms of scattering (see
\cite{ORW}).

This raises the following questions: Is there a reasonable way to
regularize (\ref{94}) with respect to heat kernels (for $s>1$)? Is
there a renormalization group prescription for understanding some
aspect of the critical limit $s\downarrow 1$? There is possibly an
important hint arising from work of Uhlenbeck (see \cite{U}). The
classical theory corresponding to (\ref{94}) is `completely
integrable', for all values of $t$ and $l$. This involves
reformulating the classical equations in terms of the connection
one-form $A=g^{-1}dg$ and a spectral parameter (which evolves into
the deformation parameter at the quantum level). In considering
the hyperfunction completion in one dimension, we abandoned
unitarity. In two dimensions it apparently is necessary to relax
unitarity and the zero-curvature condition.

\end{document}